\documentclass[12pt,a4paper]{amsart}
\usepackage{amssymb,amsxtra}
\usepackage[english,french]{babel}

\newtheorem{thm}{Theorem}[section]
\newtheorem{lem}[thm]{Lemma}

\theoremstyle{definition}

\theoremstyle{remark}

\theoremstyle{plain}

\theoremstyle{remark}

\newtheorem*{example}{Example}
\numberwithin{equation}{section}

\begin{document}

\title{ Enumeration of some particular   quadruple  Persymmetric  Matrices over $\mathbb{F}_{2} $ by rank}
\author{Jorgen~Cherly}
\address{D\'epartement de Math\'ematiques, Universit\'e de
    Brest, 29238 Brest cedex~3, France}
\email{Jorgen.Cherly@univ-brest.fr}
\email{andersen69@wanadoo.fr}

\maketitle 
\begin{abstract}
Dans cet article nous comptons le nombre de certaines  quadruples  matrices persym\' etriques de rang i sur $ \mathbb {F} _ {2} . $

 \end{abstract}

\selectlanguage{english}

\begin{abstract}
In this paper we count the number of some particular  quadruple  persymmetric rank i matrices over  $ \mathbb{F}_{2}.$
 \end{abstract}
 
  \maketitle 
\newpage
\tableofcontents
\newpage

  \section{Introduction}
  \label{sec 1}  
  In this paper we propose to compute in the most simple case  the number of  quadruple  persymmetric 
  matrices with entries in   $ \mathbb{F}_{2}$ of rank i\\
  That is to compute the number  $ \Gamma_{i}^{\left[2\atop{2 \atop{ 2\atop 2}}\right]\times k}$ 
  of  quadruple  persymmetric matrices in   $ \mathbb{F}_{2}$ of rank i  $(0\leqslant i\leqslant\inf(8,k) )$ of the below form.\\
    \begin{equation}
    \label{eq 1.1}
   \left (  \begin{array} {cccccccc}
\alpha  _{1}^{(1)} & \alpha  _{2}^{(1)}  &   \alpha_{3}^{(1)} &   \alpha_{4}^{(1)} &   \alpha_{5}^{(1)} &  \alpha_{6}^{(1)}  & \ldots  &  \alpha_{k}^{(1)} \\
\alpha  _{2}^{(1)} & \alpha  _{3}^{(1)}  &   \alpha_{4}^{(1)} &   \alpha_{5}^{(1)} &   \alpha_{6}^{(1)} &  \alpha_{7}^{(1)} & \ldots  &  \alpha_{k+1}^{(1)} \\ 
\hline \\
\alpha  _{1}^{(2)} & \alpha  _{2}^{(2)}  &   \alpha_{3}^{(2)} &   \alpha_{4}^{(2)} &   \alpha_{5}^{(2)} &  \alpha_{6}^{(2)} & \ldots   &  \alpha_{k}^{(2)} \\
\alpha  _{2}^{(2)} & \alpha  _{3}^{(2)}  &   \alpha_{4}^{(2)} &   \alpha_{5}^{(2)}&   \alpha_{6}^{(2)} &  \alpha_{7}^{(2)}  & \ldots  &  \alpha_{k+1}^{(2)} \\ 
\hline\\
\alpha  _{1}^{(3)} & \alpha  _{2}^{(3)}  &   \alpha_{3}^{(3)}  &   \alpha_{4}^{(3)} &   \alpha_{5}^{(3)} &  \alpha_{6}^{(3)} & \ldots  &  \alpha_{k}^{(3)} \\
\alpha  _{2}^{(3)} & \alpha  _{3}^{(3)}  &   \alpha_{4}^{(3)}&   \alpha_{5}^{(3)} &   \alpha_{6}^{(3)}  &  \alpha_{7}^{(3)} & \ldots  &  \alpha_{k+1}^{(3)} \\ 
\hline \\
\alpha  _{1}^{(4)} & \alpha  _{2}^{(4)}  &   \alpha_{3}^{(4)} &   \alpha_{4}^{(4)} &   \alpha_{5}^{(4)}  &  \alpha_{6}^{(4)} & \ldots  &  \alpha_{k}^{(4)} \\
\alpha  _{2}^{(4)} & \alpha  _{3}^{(4)}  &   \alpha_{4}^{(4)}&   \alpha_{5}^{(4)} &   \alpha_{6}^{(4)}  &  \alpha_{7}^{(4)} & \ldots  &  \alpha_{k+1}^{(4)} \\ 
\end{array} \right )  
\end{equation} 
We remark that this paper is based on the results in  the authors paper [11]

   \section{Some notations concerning the field of Laurent Series $ \mathbb{F}_{2}((T^{-1})) $ }
  \label{sec 2}  
  We denote by $ \mathbb{F}_{2}\big(\big({T^{-1}}\big) \big)
 = \mathbb{K} $ the completion
 of the field $\mathbb{F}_{2}(T), $  the field of  rational fonctions over the
 finite field\; $\mathbb{F}_{2}$,\; for the  infinity  valuation \;
 $ \mathfrak{v}=\mathfrak{v}_{\infty }$ \;defined by \;
 $ \mathfrak{v}\big(\frac{A}{B}\big) = degB -degA $ \;
 for each pair (A,B) of non-zero polynomials.
 Then every element non-zero t in
  $\mathbb{F}_{2}\big(\big({\frac{1}{T}}\big) \big) $
 can be expanded in a unique way in a convergent Laurent series
                              $  t = \sum_{j= -\infty }^{-\mathfrak{v}(t)}t_{j}T^j
                                 \; where\; t_{j}\in \mathbb{F}_{2}. $\\
  We associate to the infinity valuation\; $\mathfrak{v}= \mathfrak{v}_{\infty }$
   the absolute value \; $\vert \cdot \vert_{\infty} $\; defined by \;
  \begin{equation*}
  \vert t \vert_{\infty} =  \vert t \vert = 2^{-\mathfrak{v}(t)}. \\
\end{equation*}
    We denote  E the  Character of the additive locally compact group
$  \mathbb{F}_{2}\big(\big({\frac{1}{T}}\big) \big) $ defined by \\
\begin{equation*}
 E\big( \sum_{j= -\infty }^{-\mathfrak{v}(t)}t_{j}T^j\big)= \begin{cases}
 1 & \text{if      }   t_{-1}= 0, \\
  -1 & \text{if      }   t_{-1}= 1.
    \end{cases}
\end{equation*}
  We denote $\mathbb{P}$ the valuation ideal in $ \mathbb{K},$ also denoted the unit interval of  $\mathbb{K},$ i.e.
  the open ball of radius 1 about 0 or, alternatively, the set of all Laurent series 
   $$ \sum_{i\geq 1}\alpha _{i}T^{-i}\quad (\alpha _{i}\in  \mathbb{F}_{2} ) $$ and, for every rational
    integer j,  we denote by $\mathbb{P}_{j} $
     the  ideal $\left\{t \in \mathbb{K}|\; \mathfrak{v}(t) > j \right\}. $
     The sets\; $ \mathbb{P}_{j}$\; are compact subgroups  of the additive
     locally compact group \; $ \mathbb{K}. $\\
      All $ t \in \mathbb{F}_{2}\Big(\Big(\frac{1}{T}\Big)\Big) $ may be written in a unique way as
$ t = [t] + \left\{t\right\}, $ \;  $  [t] \in \mathbb{F}_{2}[T] ,
 \; \left\{t\right\}\in \mathbb{P}  ( =\mathbb{P}_{0}). $\\
 We denote by dt the Haar measure on \; $ \mathbb{K} $\; chosen so that \\
  $$ \int_{\mathbb{P}}dt = 1. $$\\
  
  $$ Let \quad
  (t_{1},t_{2},\ldots,t_{n} )
 =  \big( \sum_{j=-\infty}^{-\nu(t_{1})}\alpha _{j}^{(1)}T^{j},  \sum_{j=-\infty}^{-\nu(t_{2})}\alpha _{j}^{(2)}T^{j} ,\ldots, \sum_{j=-\infty}^{-\nu(t_{n})}\alpha _{j}^{(n)}T^{j}\big) \in  \mathbb{K}^{n}. $$ 
 We denote $\psi  $  the  Character on  $(\mathbb{K}^n, +) $ defined by \\
 \begin{align*}
  \psi \big( \sum_{j=-\infty}^{-\nu(t_{1})}\alpha _{j}^{(1)}T^{j},  \sum_{j=-\infty}^{-\nu(t_{2})}\alpha _{j}^{(2)}T^{j} ,\ldots, \sum_{j=-\infty}^{-\nu(t_{n})}\alpha _{j}^{(n)}T^{j}\big) & = E \big( \sum_{j=-\infty}^{-\nu(t_{1})}\alpha _{j}^{(1)}T^{j}\big) \cdot E\big( \sum_{j=-\infty}^{-\nu(t_{2})}\alpha _{j}^{(2)}T^{j}\big)\cdots E\big(  \sum_{j=-\infty}^{-\nu(t_{n})}\alpha _{j}^{(n)}T^{j}\big) \\
  & = 
    \begin{cases}
 1 & \text{if      }     \alpha _{-1}^{(1)} +    \alpha _{-1}^{(2)}  + \ldots +   \alpha _{-1}^{(n)}   = 0 \\
  -1 & \text{if      }    \alpha _{-1}^{(1)} +    \alpha _{-1}^{(2)}  + \ldots +   \alpha _{-1}^{(n)}   =1                                                                                                                          
    \end{cases}
  \end{align*}
   \section{Some results concerning  n-times persymmetric matrices over  $ \mathbb{F}_{2}$}
  \label{sec 3}  
     $$ Set\quad
  (t_{1},t_{2},\ldots,t_{n} )
 =  \big( \sum_{i\geq 1}\alpha _{i}^{(1)}T^{-i}, \sum_{i \geq 1}\alpha  _{i}^{(2)}T^{-i},\sum_{i \geq 1}\alpha _{i}^{(3)}T^{-i},\ldots,\sum_{i \geq 1}\alpha _{i}^{(n)}T^{-i}   \big) \in  \mathbb{P}^{n}. $$

     Denote by $D^{\left[2 \atop{\vdots \atop 2}\right]\times k}(t_{1},t_{2},\ldots,t_{n} ) $
    
    the following $2n \times k $ \;  n-times  persymmetric  matrix  over the finite field  $\mathbb{F}_{2} $ 
    
  \begin{equation}
  \label{eq 3.1}
   \left (  \begin{array} {cccccccc}
\alpha  _{1}^{(1)} & \alpha  _{2}^{(1)}  &   \alpha_{3}^{(1)} &   \alpha_{4}^{(1)} &   \alpha_{5}^{(1)} &  \alpha_{6}^{(1)}  & \ldots  &  \alpha_{k}^{(1)} \\
\alpha  _{2}^{(1)} & \alpha  _{3}^{(1)}  &   \alpha_{4}^{(1)} &   \alpha_{5}^{(1)} &   \alpha_{6}^{(1)} &  \alpha_{7}^{(1)} & \ldots  &  \alpha_{k+1}^{(1)} \\ 
\hline \\
\alpha  _{1}^{(2)} & \alpha  _{2}^{(2)}  &   \alpha_{3}^{(2)} &   \alpha_{4}^{(2)} &   \alpha_{5}^{(2)} &  \alpha_{6}^{(2)} & \ldots   &  \alpha_{k}^{(2)} \\
\alpha  _{2}^{(2)} & \alpha  _{3}^{(2)}  &   \alpha_{4}^{(2)} &   \alpha_{5}^{(2)}&   \alpha_{6}^{(2)} &  \alpha_{7}^{(2)}  & \ldots  &  \alpha_{k+1}^{(2)} \\ 
\hline\\
\alpha  _{1}^{(3)} & \alpha  _{2}^{(3)}  &   \alpha_{3}^{(3)}  &   \alpha_{4}^{(3)} &   \alpha_{5}^{(3)} &  \alpha_{6}^{(3)} & \ldots  &  \alpha_{k}^{(3)} \\
\alpha  _{2}^{(3)} & \alpha  _{3}^{(3)}  &   \alpha_{4}^{(3)}&   \alpha_{5}^{(3)} &   \alpha_{6}^{(3)}  &  \alpha_{7}^{(3)} & \ldots  &  \alpha_{k+1}^{(3)} \\ 
\hline \\
\vdots & \vdots & \vdots  & \vdots  & \vdots & \vdots  & \vdots & \vdots \\
\hline \\
\alpha  _{1}^{(n)} & \alpha  _{2}^{(n)}  &   \alpha_{3}^{(n)} &   \alpha_{4}^{(n)} &   \alpha_{5}^{(n)}  &  \alpha_{6}^{(n)} & \ldots  &  \alpha_{k}^{(n)} \\
\alpha  _{2}^{(n)} & \alpha  _{3}^{(n)}  &   \alpha_{4}^{(n)}&   \alpha_{5}^{(n)} &   \alpha_{6}^{(n)}  &  \alpha_{7}^{(n)} & \ldots  &  \alpha_{k+1}^{(n)} \\ 
\end{array} \right )  
\end{equation} 
We denote by  $ \Gamma_{i}^{\left[2\atop{\vdots \atop 2}\right]\times k}$  the number of rank i n-times persymmetric matrices over $\mathbb{F}_{2}$ of the above form :  \\

  Let $ \displaystyle  f (t_{1},t_{2},\ldots,t_{n} ) $  be the exponential sum  in $ \mathbb{P}^{n} $ defined by\\
    $(t_{1},t_{2},\ldots,t_{n} ) \displaystyle\in \mathbb{P}^{n}\longrightarrow \\
    \sum_{deg Y\leq k-1}\sum_{deg U_{1}\leq  1}E(t_{1} YU_{1})
  \sum_{deg U_{2} \leq 1}E(t _{2} YU_{2}) \ldots \sum_{deg U_{n} \leq 1} E(t _{n} YU_{n}). $\vspace{0.5 cm}\\
    Then
  $$     f_{k} (t_{1},t_{2},\ldots,t_{n} ) =
  2^{2n+k- rank\big[ D^{\left[2\atop{\vdots \atop 2}\right]\times k}(t_{1},t_{2},\ldots,t_{n} )\big] } $$

    Hence  the number denoted by $ R_{q,n}^{(k)} $ of solutions \\
  
 $(Y_1,U_{1}^{(1)},U_{2}^{(1)}, \ldots,U_{n}^{(1)}, Y_2,U_{1}^{(2)},U_{2}^{(2)}, 
\ldots,U_{n}^{(2)},\ldots  Y_q,U_{1}^{(q)},U_{2}^{(q)}, \ldots,U_{n}^{(q)}   ) \in (\mathbb{F}_{2}[T])^{(n+1)q}$ \vspace{0.5 cm}\\
 of the polynomial equations  \vspace{0.5 cm}
  \[\left\{\begin{array}{c}
 Y_{1}U_{1}^{(1)} + Y_{2}U_{1}^{(2)} + \ldots  + Y_{q}U_{1}^{(q)} = 0  \\
    Y_{1}U_{2}^{(1)} + Y_{2}U_{2}^{(2)} + \ldots  + Y_{q}U_{2}^{(q)} = 0\\
    \vdots \\
   Y_{1}U_{n}^{(1)} + Y_{2}U_{n}^{(2)} + \ldots  + Y_{q}U_{n}^{(q)} = 0 
 \end{array}\right.\]
 
    $ \Leftrightarrow
    \begin{pmatrix}
   U_{1}^{(1)} & U_{1}^{(2)} & \ldots  & U_{1}^{(q)} \\ 
      U_{2}^{(1)} & U_{2}^{(2)}  & \ldots  & U_{2}^{(q)}  \\
\vdots &   \vdots & \vdots &   \vdots   \\
U_{n}^{(1)} & U_{n}^{(2)}   & \ldots  & U_{n}^{(q)} \\
 \end{pmatrix}  \begin{pmatrix}
   Y_{1} \\
   Y_{2}\\
   \vdots \\
   Y_{q} \\
  \end{pmatrix} =   \begin{pmatrix}
  0 \\
  0 \\
  \vdots \\
  0 
  \end{pmatrix} $\\

    satisfying the degree conditions \\
                   $$  degY_i \leq k-1 ,
                   \quad degU_{j}^{(i)} \leq 1, \quad  for \quad 1\leq j\leq n  \quad 1\leq i \leq q $$ \\
  is equal to the following integral over the unit interval in $ \mathbb{K}^{n} $
    $$ \int_{\mathbb{P}^{n}} f_{k}^{q}(t_{1},t_{2},\ldots,t_{n}) dt_{1}dt _{2}\ldots dt _{n}. $$
  Observing that $ f (t_{1},t_{2},\ldots,t_{n} ) $ is constant on cosets of $ \prod_{j=1}^{n}\mathbb{P}_{k+1} $ in $ \mathbb{P}^{n} $\;
  the above integral is equal to 
  
  \begin{equation}
  \label{eq 3.2}
 2^{q(2n+k) - (k+1)n}\sum_{i = 0}^{\inf(2n,k)}
  \Gamma_{i}^{\left[2\atop{\vdots \atop 2}\right]\times k} 2^{-iq} =  R_{q,n}^{(k)} 
 \end{equation}
 
 \begin{eqnarray}
 \label{eq 3.3}
\text{ Recall that $ R_{q,n}^{(k)}$ is equal to the number of solutions of the polynomial system} \nonumber \\
    \begin{pmatrix}
   U_{1}^{(1)} & U_{1}^{(2)} & \ldots  & U_{1}^{(q)} \\ 
      U_{2}^{(1)} & U_{2}^{(2)}  & \ldots  & U_{2}^{(q)}  \\
\vdots &   \vdots & \vdots &   \vdots   \\
U_{n}^{(1)} & U_{n}^{(2)}   & \ldots  & U_{n}^{(q)} \\
 \end{pmatrix}  \begin{pmatrix}
   Y_{1} \\
   Y_{2}\\
   \vdots \\
   Y_{q} \\
  \end{pmatrix} =   \begin{pmatrix}
  0 \\
  0 \\
  \vdots \\
  0 
  \end{pmatrix} \\
 \text{ satisfying the degree conditions}\nonumber \\
                     degY_i \leq k-1 ,
                   \quad degU_{j}^{(i)} \leq 1, \quad  for \quad 1\leq j\leq n  \quad 1\leq i \leq q  \nonumber
 \end{eqnarray}

 From \eqref{eq 3.2} we obtain for q = 1\\
   \begin{align}
  \label{eq 3.4}
 2^{k-(k-1)n}\sum_{i = 0}^{\inf(2n,k)}
 \Gamma_{i}^{\left[2\atop{\vdots \atop 2}\right]\times k} 2^{-i} =  R_{1,n}^{(k)} = 2^{2n}+2^k-1
  \end{align}

We have obviously \\

   \begin{align}
  \label{eq 3.5}
 \sum_{i = 0}^{k}
 \Gamma_{i}^{\left[2\atop{\vdots \atop 2}\right]\times k}  = 2^{(k+1)n}  
 \end{align}

From  the fact that the number of rank one persymmetric  matrices over $\mathbb{F}_{2}$ is equal to three  we obtain using
 combinatorial methods  : \\
 
    \begin{equation}
  \label{eq 3.6}
 \Gamma_{1}^{\left[2\atop{\vdots \atop 2}\right]\times k}  = (2^{n}-1)\cdot 3
  \end{equation}
  For more details see Cherly  [11]
  
  \subsection{The case n=4}

       $$ Set\quad
  (t_{1},t_{2},t_{3},t_{4} )
 =  \big( \sum_{i\geq 1}\alpha _{i}^{(1)}T^{-i}, \sum_{i \geq 1}\alpha  _{i}^{(2)}T^{-i},\sum_{i \geq 1}\alpha _{i}^{(3)}T^{-i},\sum_{i \geq 1}\alpha _{i}^{(4)}T^{-i}   \big) \in  \mathbb{P}^{4}. $$

     Denote by $D^{\left[2 \atop{2\atop{2 \atop 2}}\right]\times k}(t_{1},t_{2},t_{3},t_{4}) $
    
    the following $8 \times k $ \; quadruple  persymmetric  matrix  over the finite field  $\mathbb{F}_{2} $ 
    
  \begin{displaymath}
   \left (  \begin{array} {cccccccc}
\alpha  _{1}^{(1)} & \alpha  _{2}^{(1)}  &   \alpha_{3}^{(1)} &   \alpha_{4}^{(1)} &   \alpha_{5}^{(1)} &  \alpha_{6}^{(1)}  & \ldots  &  \alpha_{k}^{(1)} \\
\alpha  _{2}^{(1)} & \alpha  _{3}^{(1)}  &   \alpha_{4}^{(1)} &   \alpha_{5}^{(1)} &   \alpha_{6}^{(1)} &  \alpha_{7}^{(1)} & \ldots  &  \alpha_{k+1}^{(1)} \\ 
\hline \\
\alpha  _{1}^{(2)} & \alpha  _{2}^{(2)}  &   \alpha_{3}^{(2)} &   \alpha_{4}^{(2)} &   \alpha_{5}^{(2)} &  \alpha_{6}^{(2)} & \ldots   &  \alpha_{k}^{(2)} \\
\alpha  _{2}^{(2)} & \alpha  _{3}^{(2)}  &   \alpha_{4}^{(2)} &   \alpha_{5}^{(2)}&   \alpha_{6}^{(2)} &  \alpha_{7}^{(2)}  & \ldots  &  \alpha_{k+1}^{(2)} \\ 
\hline\\
\alpha  _{1}^{(3)} & \alpha  _{2}^{(3)}  &   \alpha_{3}^{(3)}  &   \alpha_{4}^{(3)} &   \alpha_{5}^{(3)} &  \alpha_{6}^{(3)} & \ldots  &  \alpha_{k}^{(3)} \\
\alpha  _{2}^{(3)} & \alpha  _{3}^{(3)}  &   \alpha_{4}^{(3)}&   \alpha_{5}^{(3)} &   \alpha_{6}^{(3)}  &  \alpha_{7}^{(3)} & \ldots  &  \alpha_{k+1}^{(3)} \\ 
\hline \\
\alpha  _{1}^{(4)} & \alpha  _{2}^{(4)}  &   \alpha_{3}^{(4)} &   \alpha_{4}^{(4)} &   \alpha_{5}^{(4)}  &  \alpha_{6}^{(4)} & \ldots  &  \alpha_{k}^{(4)} \\
\alpha  _{2}^{(4)} & \alpha  _{3}^{(4)}  &   \alpha_{4}^{(4)}&   \alpha_{5}^{(4)} &   \alpha_{6}^{(4)}  &  \alpha_{7}^{(4)} & \ldots  &  \alpha_{k+1}^{(4)} \\ 
\end{array} \right )  
\end{displaymath} 
We denote by  $ \Gamma_{i}^{\left[2\atop{2 \atop{ 2\atop 2}}\right]\times k}$  the number of rank i quadruple persymmetric matrices over $\mathbb{F}_{2}$ of the above form :  \\

  Let $ \displaystyle  f (t_{1},t_{2},t_{3},t_{4} ) $  be the exponential sum  in $ \mathbb{P}^{4} $ defined by\\
    $(t_{1},t_{2},t_{3},t_{4}) \displaystyle\in \mathbb{P}^{4}\longrightarrow \\
    \sum_{deg Y\leq k-1}\sum_{deg U_{1}\leq  1}E(t_{1} YU_{1})
  \sum_{deg U_{2} \leq 1}E(t _{2} YU_{2}) \sum_{deg U_{3} \leq 1}E(t _{3} YU_{3}) \sum_{deg U_{4} \leq 1} E(t _{4} YU_{4}). $\vspace{0.5 cm}\\
    Then
  $$     f_{k} (t_{1},t_{2},t_{3},t_{4} ) =
  2^{8+k- rank\big[ D^{\left[2\atop{2\atop{2 \atop 2}}\right]\times k}(t_{1},t_{2},t_{3},t_{4} )\big] } $$

    Hence  the number denoted by $ R_{q,4}^{(k)} $ of solutions \\
  
 $(Y_1,U_{1}^{(1)},U_{2}^{(1)},U_{3}^{(1)} ,U_{4}^{(1)}, Y_2,U_{1}^{(2)},U_{2}^{(2)}, 
U_{3}^{(2)},U_{4}^{(2)},\ldots  Y_q,U_{1}^{(q)},U_{2}^{(q)}, U_{3}^{(q)},U_{4}^{(q)}   ) \in (\mathbb{F}_{2}[T])^{5q}$ \vspace{0.5 cm}\\
 of the polynomial equations  \vspace{0.5 cm}
  \[\left\{\begin{array}{c}
 Y_{1}U_{1}^{(1)} + Y_{2}U_{1}^{(2)} + \ldots  + Y_{q}U_{1}^{(q)} = 0  \\
    Y_{1}U_{2}^{(1)} + Y_{2}U_{2}^{(2)} + \ldots  + Y_{q}U_{2}^{(q)} = 0\\
    Y_{1}U_{3}^{(1)} + Y_{3}U_{3}^{(2)} + \ldots  + Y_{q}U_{3}^{(q)} = 0\\ 
   Y_{1}U_{4}^{(1)} + Y_{2}U_{4}^{(2)} + \ldots  + Y_{q}U_{4}^{(q)} = 0 
 \end{array}\right.\]
 
    $ \Leftrightarrow
    \begin{pmatrix}
   U_{1}^{(1)} & U_{1}^{(2)} & \ldots  & U_{1}^{(q)} \\ 
      U_{2}^{(1)} & U_{2}^{(2)}  & \ldots  & U_{2}^{(q)}  \\
 U_{3}^{(1)} & U_{3}^{(2)}  & \ldots  & U_{3}^{(q)}  \\
U_{4}^{(1)} & U_{4}^{(2)}   & \ldots  & U_{4}^{(q)} \\
 \end{pmatrix}  \begin{pmatrix}
   Y_{1} \\
   Y_{2}\\
   \vdots \\
   Y_{q} \\
  \end{pmatrix} =   \begin{pmatrix}
  0 \\
  0 \\
  \vdots \\
  0 
  \end{pmatrix} $\\

    satisfying the degree conditions \\
                   $$  degY_i \leq k-1 ,
                   \quad degU_{j}^{(i)} \leq 1, \quad  for \quad 1\leq j\leq 4  \quad 1\leq i \leq q $$ \\
  is equal to the following integral over the unit interval in $ \mathbb{K}^{n} $
    $$ \int_{\mathbb{P}^{4}} f_{k}^{q}(t_{1},t_{2},t_{3},t_{4}) dt_{1}dt _{2}dt_{3} dt _{4}. $$
  Observing that $ f (t_{1},t_{2},t_{3},t_{4} ) $ is constant on cosets of $ \prod_{j=1}^{4}\mathbb{P}_{k+1} $ in $ \mathbb{P}^{4} $\;
  the above integral is equal to 
  
  \begin{equation}
  \label{eq 3.7}
 2^{q(8+k) - 4(k+1)}\sum_{i = 0}^{\inf{(8,k)}}
  \Gamma_{i}^{\left[2\atop{2\atop{2 \atop 2}}\right]\times k} 2^{-iq} =  R_{q,4}^{(k)} \quad \text{where} \; k\geqslant 1
 \end{equation}

   \section{ Computation of the number of rank 4 n-times persymmetric matrices over $\mathbb{F}_{2}$ of the form \eqref{eq 3.1}}
  \label{sec 4}  
Recall (see section \ref{sec 3} ) that we denote by  $ \Gamma_{4}^{\left[2\atop{\vdots \atop 2}\right]\times k}$  the number of rank 4 n-times persymmetric matrices over $\mathbb{F}_{2}$ of the below form :  \\
   \begin{displaymath}
   \left (  \begin{array} {cccccccc}
\alpha  _{1}^{(1)} & \alpha  _{2}^{(1)}  &   \alpha_{3}^{(1)} &   \alpha_{4}^{(1)} &   \alpha_{5}^{(1)} &  \alpha_{6}^{(1)}  & \ldots  &  \alpha_{k}^{(1)} \\
\alpha  _{2}^{(1)} & \alpha  _{3}^{(1)}  &   \alpha_{4}^{(1)} &   \alpha_{5}^{(1)} &   \alpha_{6}^{(1)} &  \alpha_{7}^{(1)} & \ldots  &  \alpha_{k+1}^{(1)} \\ 
\hline \\
\alpha  _{1}^{(2)} & \alpha  _{2}^{(2)}  &   \alpha_{3}^{(2)} &   \alpha_{4}^{(2)} &   \alpha_{5}^{(2)} &  \alpha_{6}^{(2)} & \ldots   &  \alpha_{k}^{(2)} \\
\alpha  _{2}^{(2)} & \alpha  _{3}^{(2)}  &   \alpha_{4}^{(2)} &   \alpha_{5}^{(2)}&   \alpha_{6}^{(2)} &  \alpha_{7}^{(2)}  & \ldots  &  \alpha_{k+1}^{(2)} \\ 
\hline\\
\alpha  _{1}^{(3)} & \alpha  _{2}^{(3)}  &   \alpha_{3}^{(3)}  &   \alpha_{4}^{(3)} &   \alpha_{5}^{(3)} &  \alpha_{6}^{(3)} & \ldots  &  \alpha_{k}^{(3)} \\
\alpha  _{2}^{(3)} & \alpha  _{3}^{(3)}  &   \alpha_{4}^{(3)}&   \alpha_{5}^{(3)} &   \alpha_{6}^{(3)}  &  \alpha_{7}^{(3)} & \ldots  &  \alpha_{k+1}^{(3)} \\ 
\hline \\
\vdots & \vdots & \vdots  & \vdots  & \vdots & \vdots  & \vdots & \vdots \\
\hline \\
\alpha  _{1}^{(n)} & \alpha  _{2}^{(n)}  &   \alpha_{3}^{(n)} &   \alpha_{4}^{(n)} &   \alpha_{5}^{(n)}  &  \alpha_{6}^{(n)} & \ldots  &  \alpha_{k}^{(n)} \\
\alpha  _{2}^{(n)} & \alpha  _{3}^{(n)}  &   \alpha_{4}^{(n)}&   \alpha_{5}^{(n)} &   \alpha_{6}^{(n)}  &  \alpha_{7}^{(n)} & \ldots  &  \alpha_{k+1}^{(n)} \\ 
\end{array} \right )  
\end{displaymath} 

We shall need the following Lemma.\\
  \begin{lem}
\label{lem 4.1}
\begin{equation}
\label{eq 4.1}
   \Gamma_{4}^{\left[2\atop{\vdots \atop 2}\right]\times k}=   \begin{cases}
0 & \text{if  } n = 1,        \\
2^{2k+2}-3\cdot2^{k+2}+128 & \text{if   } n=2,\\
7\cdot2^{2k+2}+651\cdot2^{k+2}-22624  & \text{if   }  n = 3. 
    \end{cases}
   \end{equation} 
   \begin{equation}
  \label{eq 4.2}
     \Gamma_{4}^{\left[2\atop{\vdots \atop 2}\right]\times k}=   \begin{cases}
31\cdot 2^{4n} -45\cdot 2^{3n}-161\cdot 2^{2n+1}+51\cdot2^{n+4}-480 & \text{if   }  k=5, \\                                                                                                
31\cdot 2^{4n} +515\cdot 2^{3n}-2450 \cdot 2^{2n}+3280 \cdot2^{n}-1376   & \text{if   }  k=6. 
  \end{cases}    
 \end{equation}
\end{lem}

  \begin{proof}

  Lemma \ref{lem 4.1} follows from Daykin [3], Cherly [7]   
      
     \end{proof}  
  \begin{lem}
  \label{lem 4.2}
We postulate that :\\
\begin{align}
\label{eq 4.3}
 \displaystyle  \Gamma_{4}^{\left[2\atop{\vdots \atop 2}\right]\times k} = 31\cdot2^{4n} + \frac{35\cdot2^{k}-1210}{2}\cdot2^{3n}
+ \frac{2^{2k+2}-783\cdot2^{k}+19028}{6}\cdot 2^{2n}\\
\displaystyle +(-2^{2k+1}+269\cdot2^{k}-5744)\cdot2^n 
 +\frac{2^{2k+2}-117\cdot2^{k+2}+9440}{3}  \quad \text{for} \quad k\geqslant 5 \nonumber  
   \end{align}
\end{lem}

  \begin{proof}
  From the expressions of $ \Gamma_{4}^{\left[2\atop{\vdots \atop 2}\right]\times k} $ for k=5 and k=6  in  \eqref{eq 4.1}we assume that 
  $ \Gamma_{4}^{\left[2\atop{\vdots \atop 2}\right]\times k} $ can be written in the form :\\
$  \displaystyle  \Gamma_{4}^{\left[2\atop{\vdots \atop 2}\right]\times k} = 31\cdot2^{4n} + a(k)\cdot2^{3n}+ b(k)\cdot 2^{2n}+c(k)\cdot2^n  +d(k)  \quad \text{for} \quad k\geqslant 5$  where \\
 $\displaystyle 31 + a(k)+ b(k)+c(k) +d(k) =0 = \Gamma_{4}^{\left[2\atop{\vdots \atop 2}\right]\times k} \quad \text {for} \; n=0. $\\[0.1 cm]
 From  \eqref{eq 4.2}  for n=1,n=2 and n=3 we obtain : \\
   \begin{equation*}
 \begin{cases} 
n=0 :   \displaystyle  a(k)+ b(k)+ c(k)  +d(k) = -31 \quad \text{for} \quad k\geqslant 5\\
n=1:   \displaystyle 8\cdot a(k)+4\cdot b(k)+2\cdot c(k)  +d(k) = -496 \quad \text{for} \quad k\geqslant 5\\
n=2:    \displaystyle 64\cdot a(k)+16\cdot b(k)+4\cdot c(k)  +d(k) = 2^{2k+2}-3\cdot2^{k+4}-7808 \quad \text{for} \quad k\geqslant 5\\
n=3 :   \displaystyle 512\cdot a(k)+64\cdot b(k)+8\cdot c(k)  +d(k) = 7\cdot2^{2k+2}+651\cdot2^{k+2}-149600 \quad \text{for} \quad k\geqslant 5\\
\end{cases}
    \end{equation*}
 \begin{equation*}
   \displaystyle  \Leftrightarrow
    \begin{pmatrix}
   1 &1 & 1 & 1 \\ 
    8 & 4  & 2& 1 \\
  64 & 16 &  4& 1 \\
   512 & 64  & 8 & 1 \\
 \end{pmatrix}\displaystyle  \begin{pmatrix}
   a(k) \\
   b(k) \\
   c(k)  \\
   d(k)  \\
  \end{pmatrix} =   \begin{pmatrix}
  -31 \\
  -496 \\
  2^{2k+2}-3\cdot2^{k+4}-7808\\
   7\cdot 2^{2k+2}+651\cdot2^{k+2}-149600
  \end{pmatrix} 
    \end{equation*}
    \begin{equation*} 
\displaystyle  \Leftrightarrow
    \begin{pmatrix}
   a(k) \\
   b(k) \\
   c(k)  \\
   d(k)  \\
  \end{pmatrix} = \begin{pmatrix}
  \frac{35\cdot2^{k}-1210}{2}\\
   \frac{2^{2k+2}-783\cdot2^{k}+19028}{6} \\
  -2^{2k+1}+269\cdot2^{k}-5744\\
  \frac{2^{2k+2}-117\cdot2^{k+2}+9440}{3}  \\
  \end{pmatrix} 
  \end{equation*}
   \end{proof}

     \section{Computation of the number of quadruple  persymmetric matrices of the form \eqref{eq 1.1}  of rank i}
  \label{sec 5}  
  
We establish in the following Lemma that the $ \Gamma_{i}^{\left[2\atop{\vdots \atop 2}\right]\times k} $   where  $ 0\leq i\leq \inf(2n,k)$ (see Section \ref{sec 3})
are solutions to the below  system 
\begin{lem}
\label{lem 5.1}
 \begin{equation}
  \label{eq 5.1}
 \begin{cases} 
 \displaystyle  \Gamma_{0}^{\left[2\atop{\vdots \atop 2}\right]\times k}  = 1 \quad \text{if} \quad  k\geqslant 1 \\
\displaystyle  \Gamma_{1}^{\left[2\atop{\vdots \atop 2}\right]\times k}  = (2^{n}-1)\cdot 3 \quad \text{if} \quad  k\geqslant 2 \\
\displaystyle \Gamma_{2}^{\left[2\atop{\vdots \atop 2}\right]\times k} = 7\cdot2^{2n}+(2^{k+1}-25) \cdot 2^{n}-2^{k+1}+18 \quad \text{for} \quad k\geqslant 3\\
\displaystyle  \Gamma_{3}^{\left[2\atop{\vdots \atop 2}\right]\times k} = 15\cdot2^{3n} + (7\cdot2^k-133)\cdot2^{2n}+ (294-21\cdot 2^k) \cdot 2^{n}   -176+14\cdot2^k \quad \text{for} \quad k\geqslant 4\\
\displaystyle  \Gamma_{4}^{\left[2\atop{\vdots \atop 2}\right]\times k} = 31\cdot2^{4n} + \frac{35\cdot2^{k}-1210}{2}\cdot2^{3n}
+ \frac{2^{2k+2}-783\cdot2^{k}+19028}{6}\cdot 2^{2n}\\
\displaystyle +(-2^{2k+1}+269\cdot2^{k}-5744)\cdot2^n 
 +\frac{2^{2k+2}-117\cdot2^{k+2}+9440}{3}  \quad \text{for} \quad k\geqslant 5\\
\displaystyle  \sum_{i = 0}^{\inf(2n,k)} \Gamma_{i}^{\left[2\atop{\vdots \atop 2}\right]\times k}  = 2^{(k+1)n} \\ 
  \displaystyle  \sum_{i = 0}^{\inf(2n,k)} \Gamma_{i}^{\left[2\atop{\vdots \atop 2}\right]\times k} 2^{-i}  = 2^{n+k(n-1)}+2^{(k-1)n}-2^{(k-1)n-k}\\
  \displaystyle \sum_{i = 0}^{\inf(2n,k)} \Gamma_{i}^{\left[2\atop{\vdots \atop 2}\right]\times k} 2^{-2i}  =
   2^{n+k(n-2)}+2^{-n+k(n-2)}\cdot[3\cdot2^k-3] +2^{-2n+k(n-2)}\cdot[6\cdot2^{k-1}-6] \\
   +2^{-3n+kn}-6\cdot2^{n(k-3)-k}+8\cdot2^{-3n+k(n-2)}
\end{cases}
    \end{equation}
  \end{lem}  
  \begin{proof}
  
  Lemma \ref{lem 5.1} follows from \eqref{eq 4.3} and Cherly [11, section 5]

  \end{proof}
  \subsection{The case n=4}
  From Cherly[ 10, section 2 ] we obtain that the number of rank 8 quadruple persymmetric matrices of the  form \eqref{eq 3.1} is equal to :\\
  $ \displaystyle 2^{4}\prod_{j=1}^{4}(2^{k}-2^{8 -j}) $\\
  That is :\\
  \begin{equation}
\label{eq 5.2}
 \displaystyle  \Gamma_{8}^{\left[2\atop {2\atop {2\atop2 }}\right]\times k} = 2^{4}\prod_{j=1}^{4}(2^{k}-2^{8 -j}) 
\end{equation}

From \eqref{eq 5.1} and  \eqref{eq 5.2}  we deduce that the number of rank i quadruple persymmetric matrices over $ \mathbb{F}_{2}$ are solutions to the below system.\\
 \begin{equation}
  \label{eq 5.3}
 \begin{cases} 
 \displaystyle  \Gamma_{0}^{\left[2\atop {2\atop {2\atop2 }}\right]\times k}  = 1 \quad \text{if} \quad  k\geqslant 1 \\
\displaystyle  \Gamma_{1}^{\left[2\atop {2\atop {2\atop2 }}\right]\times k} = 45 \quad \text{if} \quad  k\geqslant 2 \\
\displaystyle \Gamma_{2}^{\left[2\atop {2\atop {2\atop2 }}\right]\times k}  =  30\cdot2^{k}+1410  \quad \text{for} \quad k\geqslant 3\\
\displaystyle \Gamma_{3}^{\left[2\atop {2\atop {2\atop2 }}\right]\times k}  = 1470\cdot 2^{k}+ 31920   \quad \text{for} \quad k\geqslant 4\\
\displaystyle  \Gamma_{4}^{\left[2\atop {2\atop {2\atop2 }}\right]\times k}  =   140\cdot 2^{2k} +42420\cdot 2^{k}+276640  \quad \text{for} \quad k\geqslant 5\\
\displaystyle  \Gamma_{8}^{\left[2\atop {2\atop {2\atop2 }}\right]\times k}  =   16\cdot 2^{4k}- 3840\cdot 2^{3k} +286720\cdot2^{2k}-7864320\cdot 2^{k} +2^{26} \quad \text{for} \quad k\geqslant 8\\
\displaystyle  \sum_{i = 0}^{\inf(8,k)} \Gamma_{i}^{\left[2\atop {2\atop {2\atop2 }}\right]\times k}= 2^{4(k+1)} \\ 
  \displaystyle  \sum_{i = 0}^{\inf(8,k)} \Gamma_{i}^{\left[2\atop {2\atop {2\atop2 }}\right]\times k} 2^{-i}  = 2^{4k-4}+255\cdot2^{3k-4}\\
  \displaystyle \sum_{i = 0}^{\inf(8,k)} \Gamma_{i}^{\left[2\atop {2\atop {2\atop2 }}\right]\times k}2^{-2i}  =2^{4k-12}+405\cdot2^{3k-11} +8085\cdot2^{2k-9}  
  \end{cases}
    \end{equation} 
 \begin{thm}  
\label{thm 5.2}  

We have whenever $k\geqslant 4:$\\
    \begin{equation}
    \label{eq 5.4}
 \Gamma_{i}^{\left[2\atop {2\atop {2\atop2 }}\right]\times k}  = \begin{cases}
1 & \text{if  } i = 0,        \\
 45 & \text{if   } i=1,\\
 30\cdot2^{k}+1410  & \text{if   }  i = 2,  \\
 1470\cdot 2^{k}+ 31920  & \text{if   }  i = 3, \\
  140\cdot 2^{2k} +42420\cdot 2^{k}+276640  & \text{if   }  i=4, \\  
  6300 \cdot 2^{2k}+630000 \cdot2^{k} -11692800 & \text{if   }  i=5, \\
 120\cdot 2^{3k}+123480\cdot2^{2k} -6142080\cdot 2^{k}  +66170880    & \text{if   }  i=6. \\
 3720\cdot 2^{3k}-416640\cdot2^{2k}+13332480\cdot2^{k} -121896960 & \text{if   }  i=7. \\
    16\cdot 2^{4k}- 3840\cdot 2^{3k} +286720\cdot2^{2k}-7864320\cdot 2^{k} +2^{26}& \text{if   }  i=8.    
  \end{cases}    
  \end{equation}
 \end{thm}
 
 \begin{proof}
 Theorem \ref{thm 5.2} follows from \eqref{eq 5.3} 
 \end{proof}
 \newpage
 
 \section{Some applications}
 \subsection{Computation of $ \Gamma_{i}^{\left[2\atop {2\atop {2\atop2 }}\right]\times k}$\; for\; $0\leqslant i \leqslant  k,\; 1\leqslant k \leqslant 8 $}
                                                                                             From \eqref{eq 5.3} and \eqref{eq 5.4} we deduce  :\\
 \begin{example}
 \textbf{The case   $ n=4,  k=1$}
 \label{eq 6.1}
  \begin{equation}
    \Gamma_{i}^{\left[2\atop {2\atop {2\atop2 }}\right]\times 1}     =  \begin{cases}
1 & \text{if  } i = 0,        \\
   255  & \text{if   } i=1.
 \end{cases}    
   \end{equation}

\textbf{The case   $ n=4, k=2$}

  \begin{equation}
   \label{eq 6.2}
         \Gamma_{i}^{\left[2\atop {2\atop {2\atop2 }}\right]\times 2}  =  \begin{cases}
1 & \text{if  } i = 0,        \\
   45 & \text{if   } i=1,\\
4050 & \text{if   }  i=2. 
    \end{cases}    
   \end{equation}

\textbf{The case   $ n=4,  k=3$}

    \begin{equation}
     \label{eq 6.3}
  \Gamma_{i}^{\left[2\atop {2\atop {2\atop2 }}\right]\times 3}   =  \begin{cases}
1 & \text{if  } i = 0,        \\
   45 & \text{if   } i=1,\\
1650  & \text{if   }  i = 2,  \\
63840 & \text{if   }  i=3. 
    \end{cases}    
   \end{equation}

 \textbf{The case   $ n=4,  k=4 $}
  \begin{equation}
   \label{eq 6.4}
   \Gamma_{i}^{\left[2\atop {2\atop {2\atop2 }}\right]\times 8}  = \begin{cases}
1 & \text{if  } i = 0,        \\
 45 & \text{if   } i=1,\\
 1890 & \text{if   }  i = 2,  \\
 55440 & \text{if   }  i = 3, \\
  991200 & \text{if   }  i=4, \\  
   \end{cases}    
  \end{equation}

 \textbf{The case   $ n=4,  k=5 $}
  \begin{equation}
   \label{eq 6.5}
   \Gamma_{i}^{\left[2\atop {2\atop {2\atop2 }}\right]\times 5}  = \begin{cases}
1 & \text{if  } i = 0,        \\
 45 & \text{if   } i=1,\\
 2370 & \text{if   }  i = 2,  \\
 78960 & \text{if   }  i = 3, \\
 1777440 & \text{if   }  i=4, \\ 
 14918400 & \text{if   }  i=5, \\
 \end{cases}    
  \end{equation}

 \textbf{The case $  n=4, k=6 $}
  \begin{equation}
   \label{eq 6.6}
   \Gamma_{i}^{\left[2\atop {2\atop {2\atop2 }}\right]\times 6}  = \begin{cases}
1 & \text{if  } i = 0,        \\
 45 & \text{if   } i=1,\\
3330 & \text{if   }  i = 2,  \\
 126000 & \text{if   }  i = 3, \\
 3564960 & \text{if   }  i=4, \\  
  54432000    & \text{if   }  i=5, \\
   210309120  & \text{if   }  i=6. \\
   \end{cases}    
  \end{equation}

 \textbf{The case   $ n=4,  k=7 $}
  \begin{equation}
   \label{eq 6.7}
   \Gamma_{i}^{\left[2\atop {2\atop {2\atop2 }}\right]\times 7}  = \begin{cases}
1 & \text{if  } i = 0,        \\
 45 & \text{if   } i=1,\\
  5250 & \text{if   }  i = 2,  \\
  220080 & \text{if   }  i = 3, \\
   8000160  & \text{if   }  i=4, \\  
  172166400  & \text{if   }  i=5, \\
    1554739200  & \text{if   }  i=6. \\
  2559836160 & \text{if   }  i=7. \\
   \end{cases}    
  \end{equation}

 \textbf{The case   $ n=4,  k=8 $}
  \begin{equation}
   \label{eq 6.8}
   \Gamma_{i}^{\left[2\atop {2\atop {2\atop2 }}\right]\times 8}  = \begin{cases}
1 & \text{if  } i = 0,        \\
 45 & \text{if   } i=1,\\
 9090 & \text{if   }  i = 2,  \\
 408240 & \text{if   }  i = 3, \\
  20311200 & \text{if   }  i=4, \\  
  562464000  & \text{if   }  i=5, \\
  8599449600   & \text{if   }  i=6. \\
 146475\cdot2^{18} & \text{if   }  i=7. \\
  315\cdot2^{26} & \text{if   }  i=8.    
  \end{cases}    
  \end{equation}
\end{example} 
 \subsection{Computation of $ R_{q,4}^{(k)} \; \text{where} \; k=1$}

 We recall that (see section  \ref{sec 3} )   $ R_{q,4}^{(1)} $denote  the number of solutions \\
  
 $(Y_1,U_{1}^{(1)},U_{2}^{(1)},U_{3}^{(1)} ,U_{4}^{(1)}, Y_2,U_{1}^{(2)},U_{2}^{(2)}, 
U_{3}^{(2)},U_{4}^{(2)},\ldots  Y_q,U_{1}^{(q)},U_{2}^{(q)}, U_{3}^{(q)},U_{4}^{(q)}   ) \in (\mathbb{F}_{2}[T])^{5q}$ \vspace{0.5 cm}\\
 of the polynomial equations  \vspace{0.5 cm}
  \[\left\{\begin{array}{c}
 Y_{1}U_{1}^{(1)} + Y_{2}U_{1}^{(2)} + \ldots  + Y_{q}U_{1}^{(q)} = 0  \\
    Y_{1}U_{2}^{(1)} + Y_{2}U_{2}^{(2)} + \ldots  + Y_{q}U_{2}^{(q)} = 0\\
    Y_{1}U_{3}^{(1)} + Y_{3}U_{3}^{(2)} + \ldots  + Y_{q}U_{3}^{(q)} = 0\\ 
   Y_{1}U_{4}^{(1)} + Y_{2}U_{4}^{(2)} + \ldots  + Y_{q}U_{4}^{(q)} = 0 
 \end{array}\right.\]
 
    $ \Leftrightarrow
    \begin{pmatrix}
   U_{1}^{(1)} & U_{1}^{(2)} & \ldots  & U_{1}^{(q)} \\ 
      U_{2}^{(1)} & U_{2}^{(2)}  & \ldots  & U_{2}^{(q)}  \\
 U_{3}^{(1)} & U_{3}^{(2)}  & \ldots  & U_{3}^{(q)}  \\
U_{4}^{(1)} & U_{4}^{(2)}   & \ldots  & U_{4}^{(q)} \\
 \end{pmatrix}  \begin{pmatrix}
   Y_{1} \\
   Y_{2}\\
   \vdots \\
   Y_{q} \\
  \end{pmatrix} =   \begin{pmatrix}
  0 \\
  0 \\
  \vdots \\
  0 
  \end{pmatrix} $\\
 satisfying the degree conditions \\
                   $$  degY_i \leq 0 ,
                   \quad degU_{j}^{(i)} \leq 1, \quad  for \quad 1\leq j\leq 4  \quad 1\leq i \leq q $$ \\
 is equal to  $$  2^{9q-8}\sum_{i = 0}^{1}  \Gamma_{i}^{\left[2\atop{2\atop{2 \atop 2}}\right]\times 1} 2^{-iq}  $$
From \eqref{eq 6.1} we then  obtain :\\

$$ R_{q,4}^{(1)} = 2^{9q-8}\sum_{i = 0}^{1}  \Gamma_{i}^{\left[2\atop{2\atop{2 \atop 2}}\right]\times 1} 2^{-iq} =  2^{9q-8}\cdot[1+255\cdot2^{-q}] =  2^{9q-8}+255\cdot 2^{8q-8}$$
 Another way to compute $ R_{q,4}^{(1)} $ is to use linear algebra to solve linear systems over  $\mathbf{F}_2. $  \\
Set  $ \displaystyle  U_{j}^{(i)} = a_{i}^{j}+b_{i}^{j}T\; and \; Y_{i}=\delta_{i} \quad  for \quad 1\leq j\leq 4  \quad 1\leq i \leq q $\\
Let M denote the following $8\times q$ matrix over   $\mathbf{F}_2. $ \\
 
           $$ \begin{pmatrix}
 a_{1}^{1} &a_{2}^{1}  & \ldots  & a_{q}^{1}  \\ 
  a_{1}^{2} &a_{2}^{2}  & \ldots  & a_{q}^{2}      \\
   a_{1}^{3} &a_{2}^{3}  & \ldots  & a_{q}^{3}   \\
a_{1}^{4} &a_{2}^{4}  & \ldots  & a_{q}^{4}  \\
 b_{1}^{1} &b_{2}^{1}  & \ldots  & b_{q}^{1}  \\ 
  b_{1}^{2} &b_{2}^{2}  & \ldots  & b_{q}^{2}      \\
   b_{1}^{3} &b_{2}^{3}  & \ldots  & b_{q}^{3}   \\
b_{1}^{4} &b_{2}^{4}  & \ldots  & b_{q}^{4}  
 \end{pmatrix} $$

Then $$ R_{q,4}^{(1)} = \sum_{l=0}^{\inf(8,q)}Card\{M \mid rank M=l\}\cdot2^{q-l} =  2^{9q-8}+255\cdot 2^{8q-8}$$
We obtain from Landsberg [1] and Fisher and Alexander [2]:\\
$$Card\{M \mid rank M=l\}=\prod_{s=0}^{l-1}\frac{(2^8-2^s)(2^q-2^s)}{2^l-2^s} \; for\; l\geqslant 1$$\\
Then $$ 2^q + \sum_{l=1}^{\inf(8,q)}\left[\prod_{s=0}^{l-1}\frac{(2^8-2^s)(2^q-2^s)}{2^l-2^s}\right]\cdot2^{q-l}=  2^{9q-8}+255\cdot 2^{8q-8}$$

\subsection{Computation of $ R_{4,4}^{(k)} \;for\;k\geqslant 1$}

 \begin{eqnarray*}
\text{ Recall that $ R_{4,4}^{(k)}$ is equal to the number of solutions of the polynomial system}\\
    \begin{pmatrix}
   U_{1}^{(1)} & U_{1}^{(2)} &  U_{1}^{(3)} & U_{1}^{(4)} \\ 
   U_{2}^{(1)} & U_{2}^{(2)}  &  U_{2}^{(3)} & U_{2}^{(4)}  \\
  U_{3}^{(1)} & U_{3}^{(2)}  &U_{3}^{(3)}  & U_{3}^{(4)}\\
U_{4}^{(1)} & U_{4}^{(2)}   &  U_{4}^{(3)} & U_{4}^{(4)} \\
 \end{pmatrix}  \begin{pmatrix}
   Y_{1} \\
   Y_{2}\\
 Y_{3}  \\
   Y_{4} \\
  \end{pmatrix} =   \begin{pmatrix}
  0 \\
  0 \\
  0 \\
  0 
  \end{pmatrix} \\
 \text{ satisfying the degree conditions}r \\
                     degY_i \leq k-1 ,
                   \quad degU_{j}^{(i)} \leq 1, \quad  for \quad 1\leq j\leq 4  \quad 1\leq i \leq 4 
 \end{eqnarray*}

From \eqref{eq 3.7} with q=4 and \eqref{eq 5.4} we obtain  whenever $ k\geqslant 3$.  :\\

   \begin{align*}
 2^{28}\sum_{i = 0}^{\inf{(8,k)}}
 \Gamma_{i}^{\left[2\atop{2\atop{2 \atop 2}}\right]\times k} 2^{-4i} =  R_{4,4}^{(k)} = \\
 2^{4k}+5400\cdot2^{3k}+3763200\cdot2^{2k}+377395200\cdot2^{k}+3674603520 
  \end{align*}
  From \eqref{eq 6.2} we get in the case k=2 : \\
   \begin{align*}
 2^{28}\sum_{i = 0}^{2}
 \Gamma_{i}^{\left[2\atop{2\atop{2 \atop 2}}\right]\times 2} 2^{-4i} =  R_{4,4}^{(2)} = 5270142976
  \end{align*}
  
  From  \eqref{eq 6.1} we get in the case k=1:  \\
   \begin{align*}
 2^{28}\sum_{i = 0}^{1}
 \Gamma_{i}^{\left[2\atop{2\atop{2 \atop 2}}\right]\times 1} 2^{-4i} =  R_{4,4}^{(1)} = 4546625536
  \end{align*}

\subsection{Computation of $ \Gamma_{5}^{\left[2\atop{\vdots \atop 2}\right]\times k} $}
We recall (see section \ref{sec 3} ) that $ \Gamma_{5}^{\left[2\atop{\vdots \atop 2}\right]\times k}$ denote the number of rank 5
n-times persymmetric matrices over  $\mathbf{F}_2 $  of the form \eqref{eq 3.1}
We shall need the following Lemma :\\

  \begin{lem}
\label{lem 6.1}
\begin{equation}
\label{eq 6.9}
   \Gamma_{5}^{\left[2\atop{\vdots \atop 2}\right]\times k}=   \begin{cases}
0 & \text{if  } n = 1,        \\
0 & \text{if  } n = 2,\\
105\cdot 2^{2k+2}-315\cdot 2^{k+5}+53760 & \text{if   } n=3,\\
6300\cdot 2^{2k}+630000\cdot 2^{k}-116928  & \text{if   }  n = 4. 
    \end{cases}
   \end{equation} 
    \begin{equation}
  \label{eq 6.10}
    \Gamma_{5}^{\left[2\atop{\vdots \atop 2}\right]\times 6}=  63 \cdot 2^{5n}-93 \cdot2^{4n} -1650 \cdot2^{3n}+5040 \cdot 2^{2n}-4128 \cdot 2^n +768   
  \end{equation}

\end{lem}
\begin{proof}
Lemma \ref{lem 6.1} follows from \eqref{eq 5.4} and Cherly[].

\end{proof}

  \begin{lem}
  \label{lem 6.2}
We postulate that :\\
\begin{align}
\label{eq 6.11}
  \displaystyle  \Gamma_{5}^{\left[2\atop{\vdots \atop 2}\right]\times k} = 63\cdot2^{5n} + (\frac{155}{4}\cdot2^{k}-2573)\cdot2^{4n}
+ (\frac{5}{2}\cdot2^{2k}-\frac{2565}{4}\cdot2^{k}+29150)\cdot2^{3n}\\
\displaystyle +\frac{1}{2}\cdot(-35\cdot2^{2k}+6265\cdot2^{k}-247520)\cdot2^{2n} 
\displaystyle +(35\cdot2^{2k}-5490\cdot2^{k}+203872)\cdot2^{n} \nonumber \\
-20\cdot2^{2k}+2960\cdot2^{k}-106752  \quad \text{for} \quad k\geqslant 6 \nonumber 
 \end{align}
\end{lem}

  \begin{proof}
  From the expression of $ \Gamma_{5}^{\left[2\atop{\vdots \atop 2}\right]\times k} $  in  \eqref{eq 6.10} for k=6  we assume that 
  $ \Gamma_{5}^{\left[2\atop{\vdots \atop 2}\right]\times k} $ can be written in the form :\\
$  \displaystyle     63\cdot2^{5n}+a(k)\cdot2^{4n} + b(k)\cdot2^{3n}+ c(k)\cdot 2^{2n}+d(k)\cdot2^n  +e(k)  \quad \text{for} \quad k\geqslant 6$ \\ where  $\displaystyle 63 + a(k)+ b(k)+c(k) +d(k) +e(k) =0 = \Gamma_{5}^{\left[2\atop{\vdots \atop 2}\right]\times k} \quad \text {for} \; n=0. $\\[0.1 cm]
  From  \eqref{eq 6.9}  for n=1,n=2,n=3 and n=4 we obtain : \\
   \begin{equation*}
 \begin{cases} 
n=0 :   \displaystyle  a(k)+ b(k)+ c(k)  +d(k)+e(k) = -63 \quad \text{for} \quad k\geqslant 6\\
n=1:   \displaystyle 2^4\cdot a(k)+2^3\cdot b(k)+2^2\cdot c(k)  +2\cdot d(k)+e(k) = -63\cdot 2^5 \quad \text{for} \quad k\geqslant 6\\
n=2:    \displaystyle 2^8\cdot a(k)+2^6\cdot b(k)+2^4\cdot c(k)  +2^2\cdot d(k) +e(k) = -63\cdot 2^{10}  \quad \text{for} \quad k\geqslant 5\\
n=3 :   \displaystyle 2^{12}\cdot a(k)+2^9\cdot b(k)+2^6\cdot c(k)  +2^3d(k) +e(k) \\
= 105\cdot 2^{2k+2}-315\cdot 2^{k+5}-2010624 \quad \text{for} \quad k\geqslant 6\\
n=4 :   \displaystyle 2^{16}\cdot a(k)+2^{12}\cdot b(k)+2^8\cdot c(k)  +2^4\cdot d(k) +e(k)\\
 = 6300\cdot 2^{2k}+630000\cdot 2^{k}-77753088 \quad \text{for} \quad k\geqslant 6\\
\end{cases}
    \end{equation*}
      \begin{equation*}
   \displaystyle  \Leftrightarrow
    \begin{pmatrix}
   1 &1 & 1 & 1&1 \\ 
    2^4 & 2^3  & 2^2&2& 1 \\
     2^8 & 2^6  & 2^4&2^2& 1 \\
      2^{12} & 2^9  & 2^6  & 2^3 & 1 \\
 2^{16} & 2^{12}  & 2^8 & 2^4 & 1 \\
 \end{pmatrix}\displaystyle  \begin{pmatrix}
   a(k) \\
   b(k) \\
   c(k)  \\
   d(k)  \\
   e(k)\\
  \end{pmatrix} =   \begin{pmatrix}
  -63 \\
  -63\cdot2^5 \\
    -63\cdot2^{10} \\
105\cdot2^{2k+2}-315\cdot2^{k+5}-2010624\\
   6300\cdot 2^{2k}+630000\cdot2^{k}-77753088
  \end{pmatrix} 
    \end{equation*}
    \begin{equation*} 
\displaystyle  \Leftrightarrow
    \begin{pmatrix}
   a(k) \\
   b(k) \\
   c(k)  \\
   d(k)  \\
   e(k) \\
  \end{pmatrix} = \begin{pmatrix}
 \frac{155}{4}\cdot2^{k}-2573\\
 \frac{5}{2}\cdot2^{2k}-\frac{2565}{4}\cdot2^{k}+29150 \\
\frac{1}{2}\cdot(-35\cdot2^{2k}+6265\cdot2^{k}-247520)\\
35\cdot2^{2k}-5490\cdot2^{k}+203872   \\
  -20\cdot2^{2k}+2960\cdot2^{k}-106752  
   \end{pmatrix} 
  \end{equation*}
    \end{proof}
  \subsection{Computation of the number  $ \Gamma_{i}^{\left[2\atop{\vdots \atop 2}\right]\times k} $   of   n-times persymmetric 
   $2n  \times k$  rank i matrices  for $7\leqslant k \leqslant 8 $}

 We establish in the following Lemma that the $ \Gamma_{i}^{\left[2\atop{\vdots \atop 2}\right]\times k} $   where  $ 0\leq i\leq \inf(2n,k)$ (see Section \ref{sec 3})
are solutions to the below  system 
\begin{lem}
\label{lem 6.3}
 \begin{equation}
  \label{eq 6.12}
 \begin{cases} 
 \displaystyle  \Gamma_{0}^{\left[2\atop{\vdots \atop 2}\right]\times k}  = 1 \quad \text{if} \quad  k\geqslant 1 \\
\displaystyle  \Gamma_{1}^{\left[2\atop{\vdots \atop 2}\right]\times k}  = (2^{n}-1)\cdot 3 \quad \text{if} \quad  k\geqslant 2 \\
\displaystyle \Gamma_{2}^{\left[2\atop{\vdots \atop 2}\right]\times k} = 7\cdot2^{2n}+(2^{k+1}-25) \cdot 2^{n}-2^{k+1}+18 \quad \text{for} \quad k\geqslant 3\\
\displaystyle  \Gamma_{3}^{\left[2\atop{\vdots \atop 2}\right]\times k} = 15\cdot2^{3n} + (7\cdot2^k-133)\cdot2^{2n}+ (294-21\cdot 2^k) \cdot 2^{n}   -176+14\cdot2^k \quad \text{for} \quad k\geqslant 4\\
\displaystyle  \Gamma_{4}^{\left[2\atop{\vdots \atop 2}\right]\times k} = 31\cdot2^{4n} + \frac{35\cdot2^{k}-1210}{2}\cdot2^{3n}
+ \frac{2^{2k+2}-783\cdot2^{k}+19028}{6}\cdot 2^{2n}\\
\displaystyle +(-2^{2k+1}+269\cdot2^{k}-5744)\cdot2^n 
 +\frac{2^{2k+2}-117\cdot2^{k+2}+9440}{3}  \quad \text{for} \quad k\geqslant 5\\
   \displaystyle  \Gamma_{5}^{\left[2\atop{\vdots \atop 2}\right]\times k} = 63\cdot2^{5n} + (\frac{155}{4}\cdot2^{k}-2573)\cdot2^{4n}
+ (\frac{5}{2}\cdot2^{2k}-\frac{2565}{4}\cdot2^{k}+29150)\cdot2^{3n}\\
\displaystyle +\frac{1}{2}\cdot(-35\cdot2^{2k}+6265\cdot2^{k}-247520)\cdot2^{2n} 
\displaystyle +(35\cdot2^{2k}-5490\cdot2^{k}+203872)\cdot2^{n}\\
-20\cdot2^{2k}+2960\cdot2^{k}-106752  \quad \text{for} \quad k\geqslant 6\\
\displaystyle  \sum_{i = 0}^{\inf(2n,k)} \Gamma_{i}^{\left[2\atop{\vdots \atop 2}\right]\times k}  = 2^{(k+1)n} \\ 
  \displaystyle  \sum_{i = 0}^{\inf(2n,k)} \Gamma_{i}^{\left[2\atop{\vdots \atop 2}\right]\times k} 2^{-i}  = 2^{n+k(n-1)}+2^{(k-1)n}-2^{(k-1)n-k}\\
  \displaystyle \sum_{i = 0}^{\inf(2n,k)} \Gamma_{i}^{\left[2\atop{\vdots \atop 2}\right]\times k} 2^{-2i}  =
   2^{n+k(n-2)}+2^{-n+k(n-2)}\cdot[3\cdot2^k-3] +2^{-2n+k(n-2)}\cdot[6\cdot2^{k-1}-6] \\
   +2^{-3n+kn}-6\cdot2^{n(k-3)-k}+8\cdot2^{-3n+k(n-2)}
\end{cases}
    \end{equation}
  \end{lem}  
  \begin{proof}
   Lemma \ref{lem 6.3} follows from \eqref{eq 6.11} and  \eqref{eq 5.1}.
    \end{proof}
We deduce from \eqref{eq 6.12} with k=7  :\\
\textbf{The case k=7}

   \begin{equation}
  \label{eq 6.13}
  \Gamma_{i}^{\left[2\atop{\vdots \atop 2}\right]\times 7}=  \begin{cases}
1 & \text{if  } i = 0,        \\
 (2^{n}-1)\cdot 3 & \text{if   } i=1,\\
7\cdot2^{2n}+231 \cdot 2^{n}-238  & \text{if   }  i = 2,  \\
15\cdot 2^{3n}+763\cdot 2^{2n}-2394\cdot 2^{n}+1616   & \text{if   }  i = 3,  \\
 31\cdot 2^{4n} +1635\cdot 2^{3n}-2610 \cdot 2^{2n}\\
 -4080 \cdot2^{n}+5024  & \text{if   }  i=4, \\  
63 \cdot 2^{5n}+2387 \cdot2^{4n} -11970 \cdot2^{3n}\\
-9520 \cdot 2^{2n}+74592 \cdot 2^n -55552 & \text{if   }  i=5, \\
    127\cdot 2^{6n}-189\cdot2^{5n} -7378\cdot 2^{4n}\\
    +24240\cdot 2^{3n}+35168 \cdot 2^{2n}-166656 \cdot2^{n} +114688 & \text{if   }  i=6. \\
   2^{8n}- 127\cdot2^{6n}+126\cdot 2^{5n}+4960\cdot 2^{4n}\\
   -13920\cdot 2^{3n}-23808 \cdot 2^{2n}+98304 \cdot2^{n}  -65536 & \text{if   }  i=7.    
  \end{cases}    
  \end{equation}
  
 We deduce from \eqref{eq 6.12} with k=8  :\\ 
  
\textbf{The case k=8}
  \begin{equation}
  \label{eq 6.14}
  \Gamma_{i}^{\left[2\atop{\vdots \atop 2}\right]\times 8} = \begin{cases}
1 & \text{if  } i = 0,        \\
 (2^{n}-1)\cdot 3 & \text{if   } i=1,\\
 7\cdot2^{2n}+487 \cdot 2^{n}-494  & \text{if   }  i = 2,  \\
15\cdot 2^{3n}+1659\cdot 2^{2n}\\
-5082\cdot 2^{n}+3408   & \text{if   }  i = 3,  \\
 31\cdot 2^{4n} +3875\cdot 2^{3n}+13454 \cdot 2^{2n}\\
-67952 \cdot2^{n}+50592  & \text{if   }  i=4, \\  
 63 \cdot 2^{5n}+7347 \cdot2^{4n} +28830 \cdot2^{3n}\\
 -468720 \cdot 2^{2n}+1092192 \cdot 2^n -659712 & \text{if   }  i=5, \\
127\cdot 2^{6n}+10227\cdot2^{5n} -52514\cdot 2^{4n}\\
-339760\cdot 2^{3n}+2548448 \cdot 2^{2n}-4804352 \cdot2^{n} +2637824& \text{if   }  i=6. \\
    255\cdot 2^{7n}-381\cdot2^{6n}-31122\cdot2^{5n} \\
    +105648\cdot 2^{4n}+758880\cdot 2^{3n}-4617984 \cdot 2^{2n}+7913472 \cdot2^{n} -4128768 & \text{if   }  i=7. \\
2^{9n}- 255\cdot 2^{7n} +254\cdot2^{6n}+20832\cdot 2^{5n}\\
-60512\cdot 2^{4n}-451840\cdot 2^{3n}+2523136 \cdot 2^{2n}-4128768 \cdot2^{n}  +2097152 & \text{if   }  i=8.    
  \end{cases}    
  \end{equation}
\subsection{Computation of the number  $  \Gamma_{i}^{\left[2\atop {2\atop {2\atop{2\atop1 }}}\right]\times k} $   of   5-times persymmetric 
   $9  \times k$  rank i matrices  for $0 \leqslant I \leqslant 9 $}
\begin{lem}
\label{lem 6.5}
The number $  \Gamma_{i}^{\left[2\atop {2\atop {2\atop{2\atop1 }}}\right]\times k}  $ of rank i matrices of the form : \\
    
  \begin{displaymath}
   \left (  \begin{array} {ccccccc}
\alpha  _{1}^{(1)} & \alpha  _{2}^{(1)}  &   \alpha_{3}^{(1)} &   \alpha_{4}^{(1)} &   \alpha_{5}^{(1)}   & \ldots  &  \alpha_{k}^{(1)} \\
\alpha  _{2}^{(1)} & \alpha  _{3}^{(1)}  &   \alpha_{4}^{(1)} &   \alpha_{5}^{(1)} &   \alpha_{6}^{(1)}  & \ldots  &  \alpha_{k+1}^{(1)} \\ 
\hline \\
\alpha  _{1}^{(2)} & \alpha  _{2}^{(2)}  &   \alpha_{3}^{(2)} &   \alpha_{4}^{(2)} &   \alpha_{5}^{(2)}  & \ldots   &  \alpha_{k}^{(2)} \\
\alpha  _{2}^{(2)} & \alpha  _{3}^{(2)}  &   \alpha_{4}^{(2)} &   \alpha_{5}^{(2)}&   \alpha_{6}^{(2)}   & \ldots  &  \alpha_{k+1}^{(2)} \\ 
\hline\\
\alpha  _{1}^{(3)} & \alpha  _{2}^{(3)}  &   \alpha_{3}^{(3)}  &   \alpha_{4}^{(3)} &   \alpha_{5}^{(3)} & \ldots  &  \alpha_{k}^{(3)} \\
\alpha  _{2}^{(3)} & \alpha  _{3}^{(3)}  &   \alpha_{4}^{(3)}&   \alpha_{5}^{(3)} &   \alpha_{6}^{(3)}  & \ldots  &  \alpha_{k+1}^{(3)} \\ 
\hline \\
\alpha  _{1}^{(4)} & \alpha  _{2}^{(4)}  &   \alpha_{3}^{(4)} &   \alpha_{4}^{(4)} &   \alpha_{5}^{(4)}   & \ldots  &  \alpha_{k}^{(4)} \\
\alpha  _{2}^{(4)} & \alpha  _{3}^{(4)}  &   \alpha_{4}^{(4)}&   \alpha_{5}^{(4)} &   \alpha_{6}^{(4)}  & \ldots  &  \alpha_{k+1}^{(4)}\\
\hline \\
\alpha  _{1}^{(5)} & \alpha  _{2}^{(5)}  &   \alpha_{3}^{(5)} &   \alpha_{4}^{(5)} &   \alpha_{5}^{(5)}  & \ldots  &  \alpha_{k}^{(5)} 
\end{array} \right )  
\end{displaymath} 
is equal to 
 \begin{equation}
 \label{eq 6.15}
 \begin{cases} 
 \displaystyle   1 \quad \text{if} \quad i=0 \\
 \displaystyle     2^k+89 \quad \text{if} \quad  i=1 \\
\displaystyle 165\cdot2^{k}+5550 \quad \text{if}\quad  i=2 \\
\displaystyle     30\cdot2^{2k}+13050\cdot2^{k}+249720 \quad \text{if}\quad i=3 \\
\displaystyle    3710\cdot2^{2k}+698880\cdot2^{k}+4170880 \quad \text{if}\quad  i=4 \\
 \displaystyle    140\cdot2^{3k}+ 241780\cdot2^{2k}+19757920\cdot2^{k}-378595840   \quad \text{if}\quad  i=5  \\
 \displaystyle     13980\cdot2^{3k}+ 8331120\cdot2^{2k}-424945920\cdot2^{k}+4609105920  \quad \text{if} \quad  i=6  \\
    \displaystyle    120\cdot2^{4k} +591960\cdot2^{3k}-67374720\cdot2^{2k}+
2165821440\cdot2^{k}-75675\cdot2^{18}  \quad \text{if} \quad  i=7 \\
 \displaystyle     7816\cdot2^{4k} -1875840\cdot2^{3k}+140062720\cdot2^{2k}-
3841720320\cdot2^{k}+977\cdot2^{25}  \quad \text{if} \quad  i=8  \\
 \displaystyle  2^4\cdot2^{5k}-7936\cdot2^{4k}+1269760\cdot2^{3k}-81264640\cdot2^{2k}+2080374784\cdot2^{k}-2^{34}  \quad \text{if} \quad i=9  
\end{cases}
 \end{equation}
whenever $ k\geqslant 4. $
 \end{lem}
\begin{proof}
From Cherly[4] Corollary 3.10 with some obviously modifications we obtain :\\
\begin{equation}
\label{eq 6.16}
  \Gamma_{i}^{\left[2\atop {2\atop {2\atop{2\atop1 }}}\right]\times k}
   = 2^{i}\cdot \Gamma_{i}^{\left[2\atop{2\atop{2 \atop 2}}\right]\times k} +(2^{k}-2^{i-1})\cdot \Gamma_{i-1}^{\left[2\atop{2\atop{2 \atop 2}}\right]\times k}
\text{for} \quad  0\leqslant i \leqslant \inf(k,9) 
\end{equation}
Combining \eqref{eq 6.16} and \eqref{eq 5.4} we get \eqref{eq 6.15}.

\end{proof}

\end{document}